\documentstyle[fleqn,twoside,plenum2]{article}

\begin{document}

\setcounter{page}{1443}

\newcount\tempcntc
\def\@citex[#1]#2{\if@filesw\immediate\write\@auxout{\string\citation{#2}}\fi
 \@tempcnta@\@tempcntb\m@ne\def\@citea{}\@cite{\@for\@citeb:=#2\do
  {\@ifundefined
   {b@\@citeb}{\@citeo\@tempcntb\m@ne\@citea\def\@citea{,}{\bf ?}\@warning
   {Citation `\@citeb' on page \thepage \space undefined}}%
  {\setbox\z@\hbox{\global\tempcntc0\csname b@\@citeb\endcsname\relax}%
   \ifnum\tempcntc=\z@ \@citeo\@tempcntb\m@ne
    \@citea\def\@citea{,}\hbox{\csname b@\@citeb\endcsname}%
   \else
    \advance\@tempcntb\@ne
    \ifnum\@tempcntb=\tempcntc
    \else\advance\@tempcntb\m@ne\@citeo
    \@tempcnta\tempcntc\@tempcntb\tempcntc\fi\fi}}\@citeo}{#1}}
\def\@citeo{\ifnum\@tempcnta>\@tempcntb\else\@citea\def\@citea{,\ }
 \ifnum\@tempcnta=\@tempcntb\the\@tempcnta\else
  {\advance\@tempcnta\@ne\ifnum\@tempcnta=\@tempcntb \else \def\@citea{--}\fi
   \advance\@tempcnta\m@ne\the\@tempcnta\@citea\the\@tempcntb}\fi\fi}

\def\thebibliography{%
  \tthebibliography}
\let\endthebibliography=\endlist
\def\tthebibliography#1{\normalsize
  \list{\biblabel{\arabic{enumi}}}{\settowidth\labelwidth{\biblabel{#1}}
    \labelsep 14pt \itemindent 0pt
    \leftmargin\labelwidth \advance\leftmargin\labelsep
    \itemsep  0.0\baselineskip plus 0.1\baselineskip minus 0.1\baselineskip 
    \usecounter{enumi}\let\p@\empty
    \def\theenumi{\arabic{enumi}}}%
    \def\newblock{\hskip 0.11em plus 0.33em minus -0.07em}
     \hbadness10000 \sfcode`\.=1000\relax}

\def\biblabel#1{\hfill #1.}

\def\ps{$\phantom{^*}+^*$}
\def\pb{$\phantom{^{**}}+^{**}$}
\def\secondstatemento#1#2{{\bf #1 }{\sl #2}\smallskip}
\def\statemento#1#2{\smallskip{\bf #1 }{\sl #2}\smallskip}
\def\statement#1#2{\medskip{\bf #1}~~{\sl #2}\medskip}

\renewcommand{\TH}[2]{\medskip{\bf THEOREM~${\bf #1}$.}~~{\sl #2}\medskip}
\renewcommand{\LE}[2]{\medskip{\bf LEMMA~${\bf #1}$.}~~{\sl #2}\medskip}
\renewcommand{\PRP}[2]{\smallskip{\bf Proposition~${\bf #1}$.}~~%
{\sl #2}\smallskip}
\renewcommand{\CRL}[2]{\medskip{\bf COROLLARY~#1.}~~{\sl #2}\medskip}
\newcommand{\PCO}[1]{{\bf Proof of Corollary~#1.~}} 

\def\D{\Delta}
\def\G{\Gamma}
\def\la{\lambda}
\def\f{\varphi}
\def\e{\varepsilon}
\def\l{\ell}
\def\FF{\mathop{\cal F}\nolimits}
\def\GG{\mathop{\cal G}\nolimits}
\def\TT{\mathop{\cal T}\nolimits}
\def\NN{\mathop{\cal N}\nolimits}
\def\XX{\mathop{\cal X}\nolimits}
\def\eql    {\mathop{=}      \limits}
\def\ul     {\mathop{\rm u}  \limits}
\def\ulo{\ul^{\scriptscriptstyle \circ}}
\def\capl   {\mathop{\cap}   \limits}
\def\cupl   {\mathop{\cup}   \limits}
\def\maxl   {\mathop{\max}   \limits}
\def\suml   {\mathop{\sum}   \limits}
\def\liml   {\mathop{\lim}   \limits}
\def\bigcupl{\mathop{\bigcup}\limits}
\newsymbol\smallsetminus 2072
\def\setminus#1#2{#1\kern-.08ex\smallsetminus\kern-.08ex#2}
\def\hdots{\ldots}
\def\1n{1,\hdots,n}
\def\_#1{\mathop{\hspace{-2pt}^{}_{#1}}}
\def\epr{\hfill$\square$\medskip\par}
\def\S{\Sigma}
\def\cdc{,\ldots,}
\def\parr{\newline\indent}
\def\R{\mathop{\Bbb R}{}}
\def\oj {\mathop {\bar{J}}\nolimits}
\def\Ap {\mathop {A^{\scriptscriptstyle +}}}
\def\Lp {\mathop {L^{\scriptscriptstyle +}}}
\def\lpij {\mathop {\ell_{ij}^{\scriptscriptstyle +}}}
\def\intercal{\mathop {\scriptscriptstyle T}\nolimits}

\Btit{}

\tit{519.173:512.643.8}
    {P.~Yu.~Chebotarev and E.~V.~Shamis}
    {On proximity measures for graph vertices}
    {\footnote[1]{This work was supported by the Russian Foundation
                  for Basic Research, Grant No.~96-01-01010. Partial research
                  support from the European Community under Grant
                  INTAS-96-106 is also gratefully acknowledged.}}

\translation{Institute of Control Sciences, Russian Academy of
Sciences, Moscow. Translated from Avtomatika i Telemekhanika,
No.~10, pp.~113--133, October, 1998.
Original article submitted December 24, 1997.}

\abst{
We study the properties of several proximity measures for the vertices of
weighted multigraphs and multidigraphs. Unlike the classical distance for the
vertices of connected graphs, these proximity measures are applicable to
weighted structures and take into account not only the shortest, but
also all other connections, which is desirable in many applications. To
apply these proximity measures to unweighted structures, every edge
should be assigned the same weight which determines the proportion of
taking account of two routes, from which one is one edge longer than
the other. A topological interpretation is obtained for the
Moore--Penrose generalized inverse of the Laplacian matrix of a
weighted multigraph.  }

\par\vskip24truept

\arcsect{1}{INTRODUCTION}

Proximity measures for the vertices of directed and undirected graphs
arise in many applied settings. The range of applications of such functions
is rather wide, including chemistry
%
[1--7],
crystallography
\cite{Kas}, epidemiology \cite{Altm}, urban planning \cite{Hugh},
organizational management \cite{Ibar}, political sciences \cite{HeGl},
aggregation of preferences \cite{Sha,ChSh}, etc. The most steadfast interest
in them is displayed in mathematical sociology
%
[15--25] in connection with the problem of measuring centrality in social networks. This important concept is
multifarious, and a great variety of model and heuristic approaches were proposed to define its numerical
representation. Note that graph theorists mainly dealt with the classical distance between the vertices of a
connected graph \cite{BucHar}, which is the length of the shortest path between them. At the same time, the
presence of additional, even longer paths is of practical importance in many applications. For example, if
the shortest road between two places is congested, a portion of goods can be delivered by a longer path
(detour).

In this paper, we study the properties of several ``sensitive'' proximity measures that take into account all
connections in a multigraph. Their common feature is the calculation (with appropriate weights) of all
structures of a certain type that connect two vertices: paths, routes, routes with drains, trees, and so
forth. For these measures, the weights of edges determine the proportion of taking account of longer paths in
comparison with shorter ones. In some cases, the weight of an edge has the meaning of a ``transfer factor''
that specifies the losses (of substance, influence, reliability, etc.) when moving through a graph.

\arcsect{2}{SOME NORMATIVE PROPERTIES OF PROXIMITY MEASURES}
\label{prop}

Suppose that $G$ is a weighted multigraph with vertex set $V(G)=\{\1n\}$ and edge set $E(G)$; $\G$ is a
weighted multidigraph with vertex set $V(\G)=\{\1n\}$ and arc set $E(\G)$; the weights of edges and arcs are
denoted by $\e_{ij}^p$ (the $p$th edge/arc from $i$ to $j$) and are strictly positive. The terms ``graph''
and ``subgraph'' will be used as generic ones (allowing multiple, weighted, and directed arcs).

Suppose that $E=(\e\_{ij})$ is the matrix of total weights of edges
(arcs) for all pairs of vertices:
\[
\e\_{ij}=\suml_{p=1}^{a\_{ij}}\e_{ij}^p,\;\:i,j=\1n,
\]
where $a\_{ij}$ is the number of edges (arcs) that connect $i$ to $j$. Let
$H$ be a subgraph of $G$. The product of the weights of all edges of $H$ will
be termed the {\em weight\/} of $H$ and denoted by $\e(H)$. The weight of a
directed subgraph of $\G$ is defined similarly. The weight of a subgraph
without edges/arcs is set to be~1. For any nonempty set of subgraphs $\GG$,
its weight is
\begin{equation}
\e(\GG)=\suml_{H\in\GG}\e(H).
\label{eset}
\end{equation}
The weight of the empty set is zero. $P=(p\_{ij})$ will designate various
$n\times n$-matrices of proximity (accessibility, connectedness) measures for
the vertices of $G$ or $\G$.

Let us formulate a number of conditions whose fulfillment is rather natural to
the proximity measures under consideration. Most of them were introduce in
connection with the {\it relative forest accessibility\/} of graph
vertices~\cite{weAiT97}.

\statemento{Symmetry.}{For every multigraph$,$ the matrix $P$ is symmetric.}

This condition is hardly natural as applied to directed graphs. In the
statements below, symmetry always stands for that applied in the undirected
case.

\statemento{Nonnegativity.}{For any multigraph $($multidigraph$),$
$p\_{ij}\ge0,\;\:i,j=\1n$.}

\secondstatemento{Reversal property.}{For any multidigraph, the
reversal of all its arcs $($provided that their weights are preserved$)$
results in the transposition of the proximity matrix.}

\secondstatemento{Diagonal maximality.}{For any multigraph $($multidigraph$)$
and any $i,j=\1n$ such that $i\ne j$, $\;p\_{ii}>p\_{ij}$ and
$p\_{ii}>p\_{ji}$ hold.}

This condition requires a stronger relation of each vertex to itself than to
any other vertex. If a proximity measure has the reversal property, then
the two inequalities of the diagonal maximality are equivalent in the
case of directed graphs as well as in the undirected case. Since all
the measures applicable to directed graphs hereinafter possess the
reversal property, we will prove only the first inequality of diagonal
maximality.

\statemento{Triangle inequality for proximities.}{For any multigraph
and for any $i,j,k=\1n$, $\;p\_{ij}+p\_{ik}-p\_{jk} \le p\_{ii}$ holds.
If, in addition, $j=k$ and $i\ne j$, then the inequality is strict.}

The triangle inequality for proximities is also meaningful as applied to directed graphs.  However, in this
case it requires  special consideration, since different orders of subscripts ($p\_{ij}$ or $p\_{ji}$, etc.)
give rise to several modifications. In this paper, a ``directed'' triangle inequality for proximities is used
in some proofs, but in the main text we deal with its undirected version only.

Consider the index%
\footnote{Transformations of the form of (\ref{metr}) in either explicit or
implicit form appear in many papers, e.g.,
\cite{DrbGol1,KleRan,Kle97,Altm,SteZel,weAiT97,Fiedler95,Merr97}, and also in
the theory of linear statistical models.}
\begin{equation}
d\_{ij}=p\_{ii}+p\_{jj}-p\_{ij}-p\_{ji},\quad i,j=\1n.
\label{metr}
\end{equation}

\secondstatemento{Metric representability of proximity.}{The index $d\_{ij}$
is a distance between the vertices of a multigraph, i.e., it satisfies the
axioms of a metric.
}

This condition is always satisfied, provided that symmetry and
the triangle inequality for proximities hold true \cite{weAiT98}; the
latter condition turns out to be closely related to the usual triangle
inequality for the distance $d\_{ij}$. Moreover, some kind of duality has been
established between the metrics defined on an arbitrary set and the functions
that satisfy the triangle inequality for proximities and an additional
normalization condition~\cite{weAiT98}.

Let us adduce an example not dealing with graphs to illustrate the triangle
inequality for proximities and the metric (\ref{metr}). Let $p(x,y)$ be the
function, defined on the pull-back of some family $\XX$ of finite sets, that
takes every pair of sets $(x,y)$ to the number $|x\cap y|$ of elements in
their meet. Then, for any $x,y,z\in\XX$,
\begin{eqnarray}
p(x,x)&=&|x|\ge |x\cap y|+|x\cap z|-|x\cap y\cap z|
            \ge |x\cap y|+|x\cap z|-|y\cap z|
\nonumber\\
      &=&p(x,y)+p(x,z)-p(y,z),
\label{caps}
\end{eqnarray}
i.e., the triangle inequality for proximities is fulfilled (since the first
inequality in (\ref{caps}) is strict at $x\ne y$ and $y=z$).
The transformation (\ref{metr}) applied to $p(x,y)$ generates the usual
metric on finite sets: the distance between $x$ and $y$ is the number of
elements in their symmetric difference.

In the sequel, we assume that there is one path of length $0$ from any
vertex to itself.

\statemento{Disconnection condition.}{For any multigraph $G$ $($multidigraph
$\G)$ and for any $i,j=\1n,\;$ $p\_{ij}=0$ iff there is no path from $i$
to~$j$ in $G$ $($in $\G)$.
}

\secondstatemento{Connectivity condition}{$($a consequence of the
disconnection condition$).$
\parr{\rm (1)} For any multigraph$,$ the matrix $P$ can be reduced to a
block-diagonal form$,$ where all block entries are strictly positive$,$ all
other entries being zero. The matrix $P$ is strictly positive iff $G$ is
connected.
\parr{\rm (2)} For any $i,j,k\in V(G)$, $p\_{ij}>0$ and $p\_{jk}>0$ imply
$p\_{ik}>0$.
}

The following normative property can be considered as an extension of diagonal
maximality.

\statemento{Transit property.}
{For any multigraph $G$ and any $i,k,t\in V(G),$ if $G$ contains a
path from $i$ to $k,$ $i\ne k\ne t,$ and each path from $i$ to $t$ includes
$k,$ then $p\_{ik}>p\_{it}.$ The same applies to multidigraphs.
}

\secondstatemento{Monotonicity.}{
Suppose that the weight of some edge $($arc$)$ $\e_{kt}^p$ in a multigraph $G$
{$($}multidigraph $\G)$ increases or a new edge $($arc$)$
from $k$ to $t$ appears. Then
\parr
{\rm (1)} $\D p\_{kt}>0,$ and for any $i,j=\1n,$  $\{i,j\}\ne\{k,t\}$ implies
$\D p\_{kt}>\D p\_{ij};$ in the directed case, the hypothesis is weakened to
$[i\ne k$ or $j\ne t];$
\parr
{\rm (2)} for any $i=\1n,$ if there is a path from $i$ to $k,$ and each path
from $i$ to $t$ includes $k,$ then $\D p\_{it}>\D p\_{ik};$
\parr
{\rm (3)} for any $i\_1,i\_2=\1n,$ if $i\_1$ and $i\_2$ can be substituted
for $i$ in the hypothesis of item~$2,$ then $p\_{i\_1\,i\_2}$ does not
increase.  }

Item~3 can be interpreted as follows: the proximity between two vertices
does not increase whenever the bond that appears or becomes stronger is
extraneous for the connection of these two vertices.

\arcsect{3}{PATH ACCESSIBILITY}
\label{SecPath}

The simplest proximity measure that takes into account not only the shortest path between vertices is {\it
path accessibility}. The path accessibility of $j$ from $i$ is defined as the total weight of all paths from
$i$ to $j$. There are two ways of defining this measure at $j=i$. First, ``paths from $i$ to $i$'' can be
interpreted as simple cycles from $i$ to $i$ plus the path of length $0$ whose weight is unity. The second
possibility is to assume that the latter trivial path is the only path from $i$ to $i$. Note that discarding
this trivial path leaves no chance of meeting diagonal maximality. We adopt the first definition, which is
more informative, though more disputable, but the subsequent discussion is applicable to the second
definition too.

Path accessibility can serve as a proximity measure only if a shorter path is assigned a greater weight than
a covering longer path (cf.\ transit property). If the weight of a path is the product of the weights of the
constituent edges/arcs (as we assume hereinafter), this requires that the edge/arc weights belong to the
interval $[0,1]$. In this way, path accessibility (as well as the subsequent indices) corresponds to the
models where every edge weight is a ``transfer factor'' that determines the weakening of ``vertex influence''
with movement away from the vertex along the edge. In some cases, such a model can be applicable to
transformed data that result after multiplying each edge (arc) weight by a constant factor $\tau$,
$0<\tau<(\max_{i,j,p}\e_{ij}^p)^{-1}$. With the same effect, the weight of a path can be defined as
$\prod(\tau\e(e))$, with the product over all edges (arcs) $e$ in the path. In the same manner, each edge/arc
of an unweighted graph can be assigned the same weight $\tau$.  While talking about edge/arc weights, we will
have in mind the weights so obtained too.

To choose $\tau$ for unweighted graphs, one has to estimate the proximity of two vertices connected by an
edge compared to the proximity of two vertices connected by a two-edge path. If the latter vertices appear to
be two times  ``farther,'' then $\tau=1/2$ can be chosen. In this case, two vertices connected by a
three-edge (four-edge) path are four times (respectively, eight times) farther. If the respective decrements
of 3 and 4 seem to be more natural, one has to take another model, which the reader can easily construct.
Here, the reciprocal weight of a path is the sum of the reciprocal weights of the constituent edges (harmonic
rather than geometric decrease). The original concept in such models is distance, whereas proximity can be
introduced as the reciprocal value.  Undoubtedly, these models are natural, but we do not consider them in
this paper.  Some of their properties are discussed in \cite{Altm,KirNeuSha}.

Let $P$ be the matrix whose entries are the values of path accessibility for
all pairs of vertices.

\PRP{1}{
Path accessibility has the following properties$:$
symmetry$,$
nonnegativity$,$
reversal property$,$ and
disconnection condition.
Moreover$,$ if $\e_{ij}^p<\e\_0$ for all $i,j=\1n,$ $p\le a\_{ij}$
$($where $\e\_0$ is a specific constant dependent on $n$ and the
greatest possible number $m$ of multiple edges/arcs$),$ then
diagonal maximality$,$
the triangle inequality for proximities$,$
the transit property$,$ and
monotonicity
are true.
}

The proofs of all statements are given in the Appendix.

Since the changes in proximities under special modifications of graphs
are of interest, restrictions on the edge/arc weights are
introduced for certain families of graphs rather than for individual
graphs. Here, such a family is determined by $n$ and $m$.

\arcsect{4}{CONNECTION RELIABILITY AS A VERTEX PROXIMITY MEASURE}

Let us assume that all edge/arc weights belong to the interval $[0,1]$, and consider them as the
probabilities of edge/arc intactness. Define $p\_{ij}$ to be the reliability of connections between $i$ and
$j$, i.e., the probability that at least one intact path between $i$ and $j$ survives, provided that all
edge/arc failures are independent; let $P=(p\_{ij})$ be the matrix of {\em connection reliabilities\/} for
all pairs of vertices. Connection reliability can be considered as a proximity measure for graph vertices.
Let us point out some advantages of this measure. First, it is based upon a natural model. Second, it is not
always appropriate that the proximity be doubled as all paths between a pair of vertices are duplicated (this
is the case when path accessibility is used); in some cases, the increase should be more moderate. This
property features connection reliability.

According to a well-known theorem (see, e.g., \cite[p.~10]{Shier}),
\begin{equation}
\!\!\!
p\_{ij}(G)=\suml_k \Pr(R_k)
          -\suml_{k<t}\Pr(R_k R_t)
          +\suml_{k<t<l}\Pr(R_k R_t R_l)-\ldots
          +(-1)^{h+1}\Pr(R_1 R_2\cdots R_h),
\label{P_ways}
\end{equation}
where $R_1, R_2\cdc R_h$ are all paths between $i$ and $j$;
$\Pr(R_kR_t)=\e(R_k\cup R_t)$, where $R_k\cup R_t$ is the subgraph that
contains those edges (arcs) that belong to $R_k$ or $R_t$, and so forth.
By virtue of (\ref{P_ways}), connection reliability is
a natural modification of path accessibility that takes into account the
{\it degree of overlapping\/} for different paths between two vertices.

Connection reliability possesses all the normative properties
listed in Sec.~2,
though for some of them the strict inequality $\e_{ij}^p<1$
is necessary.

\PRP{2}{
Connection reliability has the following properties$:$
symmetry$,$
nonnegativity$,$
reversal property$,$
disconnection condition$,$ and
item~$3$ of monotonicity.
Diagonal maximality$,$
the triangle inequality for proximities$,$
the transit property$,$ and
items~$1$ and $2$ of monotonicity
hold true, provided that the intactness probability of each edge/arc is
strictly less than~$1;$ otherwise they are satisfied in a nonstrict form.
}

\arcsect{5}{ROUTE ACCESSIBILITY}

A special feature of path accessibility (which also applies to connection reliability) is the necessity of a
logical algorithm for its calculation. The replacement of paths by routes reduces the problem to the
inversion of a matrix (see, e.g.,~\cite{Kas}). Moreover, the route accessibility of $j$ from $i$ has some
relation to the following problem: find the probability that a random walk started at $i$ is located at $j$
at a ``randomly chosen'' moment. Note that the proximity measures originating from the analysis of Markov
chains require special consideration. Interesting information on them can be found in
\cite{Kle97,Altm,Frie,DoyleSnell}.

Consider the matrix $P=(I-E)^{-1}$, where $E=(\e\_{ij})$ is the matrix of
total weights of edges (arcs) introduced above. Expand $P$ as the sum of an
infinitely decreasing geometric progression (not specifying the conditions
of its validity so far):
\begin{equation}
P=(I-E)^{-1}=I+E+E^2+\ldots\quad.
\label{expand}
\end{equation}
Let $\NN_{ij}$ be the set of routes from $i$ to $j$. Since the entries
of $E^k$ are the total weights of $k$-length routes, (\ref{expand})
implies \begin{equation} p\_{ij}=\suml_{N\in\NN_{ij}}\e(N),
\label{P_entry}
\end{equation}
i.e., $p\_{ij}$ is the total weight of routes from $i$ to $j$ (at $j=i,$
the route of length $0$ weighted by $1$ is naturally taken into account).
Therefore, $P$ is the matrix of {\em route accessibilities\/} in a multigraph
(multidigraph).

Equation~(\ref{expand}) is valid if and only if
\begin{equation}
\vert\la\_1\vert<1,
\label{lamb}
\end{equation}
where $\vert\la\_1\vert$ is the spectral radius of $E$
\cite[Corollary 5.6.16]{HoJo}.

Consider the upper bound for $\vert\la\_1\vert$ provided by
the Ger\v{s}gorin theorem (see \cite{HoJo}):
\begin{equation}
\vert\la\_1\vert
\le\maxl_{1\le i\le n}\suml_{j=1}^n\vert\e\_{ij}\vert.
\label{Gers}
\end{equation}

Let $\e\_{\max}$ be an imposed upper bound for the edge/arc weights; suppose
that $m$ is the greatest possible number of multiple edges (arcs) incident to
the same pair of vertices. Then
\begin{equation}
\vert\la\_1\vert
\le\maxl_{1\le i\le n}\suml_{j=1}^n\vert\e\_{ij}\vert \le m(n-1)\e\_{\max}.
\end{equation}
Therefore, the validity of (\ref{lamb}) (and thus of (\ref{expand})) is
provided by
\begin{equation}
\e\_{\max}<\big (m (n-1) \big)^{-1}.
\label{e_max}
\end{equation}

While on the subject of route accessibility, we will assume that the
constraint (\ref{e_max}) is satisfied (possibly, after the transformation of
edge/arc weights mentioned in Sec.~3).
A representation of
the entries of $P$ through the weights of specific connections in a digraph
(this representation involves finite sums only and thus does not require any
restrictions on the edge/arc weights) can be found in~\cite{Ponst}. A useful
review of results related to the calculation of routes in graphs is given in
~\cite{CveDoZa}.

\PRP{3}{
Route accessibility has the following properties$:$
symmetry$,$
nonnegativity$,$
reversal property$,$
diagonal maximality$,$
the triangle inequality for proximities
$($for the edge/arc weights not exceeding $(mn)^{-1}),$
the disconnection condition$,$
the transit property$,$ and
items~$1$ and~$2$ of monotonicity$.$
Item~$3$ of monotonicity is not valid for it.
}

The triangle inequality has not yet been proved in the general case.
The following proposition is used in the proofs of other properties and is
worth mentioning in itself.

\statemento{Proposition~4}
{
{\rm (on one-step increment of route accessibility for multidigraphs).}
Suppose that some arc weight $\e_{kt}^p$ in $\G$ increases by $\D\e\_{kt}>0$
or an extra arc from $k$ to $t$ with a weight $\D\e\_{kt}$ is added to $\G.$
Let $\G'$ be the new multidigraph and $P'=P(\G')$. Then
\[
\D P=hR,
\]
where $\D P=P'-P,$ $h=\frac{\textstyle{\D\e\_{kt}}}
                           {\textstyle{1-\D\e\_{kt} p\_{tk}}}$,
and
$R=(r\_{ij}) $ is the $n\times n$-matrix with entries
$r\_{ij}=p\_{ik}p\_{tj}.$
}

\arcsect{6}{RELATIVE FOREST ACCESSIBILITY FOR MULTIGRAPHS}

The notion of relative forest accessibility for multigraphs and multidigraphs
was introduced in \cite{weAiT97,ChShAtl}, where we studied its properties
in the case of multigraphs. In the present paper, we consider the undirected
case too. Relative forest accessibility for multidigraphs is not one,
but two complementary indices, calculated by counting the weights of
converging and diverging spanning forests, respectively. None of the two
possesses the reversal property of Sec.~2,
but they have it ``together'': the matrix of the first index for the multidigraph with reversed arcs equals
the transposed matrix of the second index for the original multidigraph, and vice versa. Some other
properties are also natural to apply to the pair of indices. Thereby, the consideration of the
above-mentioned indices in this paper could excessively complicate its structure. In the next two sections,
we study the limit properties of the relative forest accessibility measure for multigraphs. The corresponding
limit properties for multidigraphs are substantially different, and they should be considered elsewhere.

All assertions of Proposition~5 stated below, except for item~1 of
monotonicity, are proved in~\cite{weAiT97}. Item~1 of monotonicity is proved
in the Appendix.

Recall that the {\em Laplacian matrix\/} (also called the Kirchhoff or the
admittance matrix) of a multigraph $G$ is the $n\times n$-matrix
$L=L(G)=\left(\l_{ij}\right)$ with entries
\begin{eqnarray}
\label{lij}
\l_{ij}&=&-\suml_{p=1}^{a\_{ij}}\e\_{ij}^p,\;\:j\ne i,\;\:i,j=\1n,
\\
\l_{ii}&=&-\suml_{j\ne i}\l_{ij},\;\:i=\1n,
\label{lii}
\end{eqnarray}
where $a\_{ij}$ is the number of (multiple) edges incident to $i$ and
$j$ simultaneously. By (\ref{lij}) and (\ref{lii}), $\l_{ii}$ is the total
weight of edges incident to $i$ (exclusive of loops).

The matrix
\[
Q=(q\_{ij})=(I+L(G))^{-1}.
\]
is the matrix of {\em relative forest accessibilities\/} of vertices in $G$.

This term is suggested by the matrix-forests theorem
\cite{Sha,weAiT97,Merr97,ChShAtl}. Suppose that $\FF(G)=\FF$ is the set of all
spanning rooted forests of multigraph $G$, and $\FF^{ij}(G)=\FF^{ij}$ is the
set of those spanning rooted forests, in which $i$ and $j$ belong to the same
tree rooted at $i$. A {\em spanning rooted forest\/} is an acyclic subgraph of
$G$ that has the same vertex set as $G$ and one marked vertex (a root)
in each component.

\statement{THEOREM~1~(matrix-forest theorem for weighted multigraphs)}
{{\rm \cite{weAiT97,ChShAtl}}.
For any weighted multigraph $G$, the matrix $Q=(I+L(G))^{-1}$ exists
and $q\_{ij}=\e(\FF^{ij})\big/\e(\FF),\;$ $i,j=\1n$.
}

Recall that, according to (\ref{eset}), $\e(\FF^{ij})$ and $\e(\FF)$
are the total weights of forests that belong to $\FF^{ij}$ and $\FF$,
respectively.  For the sake of unification, in the sequel we denote the
matrix $Q$ by $P=(p\_{ij})$ (as well as other matrices of proximity
measures).

The characteristic features of relative forest accessibility are {\em doubly
stochastic normalization\/} (more precisely, its second condition) and {\em
macrovertex independence}.

\statemento{Doubly stochastic normalization.}{
For any multigraph $G$,
\parr{\rm (1)} $p\_{ij}\ge0,\;\:i,j=\1n$, and
\parr{\rm (2)} $\suml_{i=1}^n p\_{ij}=\suml_{i=1}^n p\_{ji}=1,\;\:i=\1n$.}

According to this condition, $p\_{ij}$ can be interpreted as the share of the connectivity of $i$ and $j$ in
the total connectivity of $i$ (or $j$) with all vertices. This interpretation requires some explanation.
Indeed, by virtue of symmetry, it requires that the ``total connectivity'' of all vertices be identical,
irrespective of the difference in their position within a multigraph. This is realized with the aid of the
diagonal entries of the matrix: if $i$ is poorly connected with other vertices, then $p\_{ii}$ (which
expresses the ``solitariness'' of $i$) is great, and hereby the ``total connectivity'' is the same as for all
other vertices.

Let $D$ be a subset of the vertex set $V(G)$. We say that $D$ is a {\em
macrovertex\/} in $G$, if for every $i,j\in D$ and $k\notin D$,
$\;\e\_{ik}=\e\_{jk}$ holds.

The following property is a sufficient condition for the equality and
stability of proximities.

\statemento{Macrovertex independence.}{Suppose that $D$ is a macrovertex in
$G$ and $i\in D,$ $j\in D,$ $k\notin D.$ Then $p\_{ik}=p\_{jk},$ and
$\:p\_{ik}$ does not vary when any new edges appear or the weights of any
existing edges change inside $D.$
}

Macrovertex independence substantially strengthens the following simple
condition (which is not included in the list of Sec.~2,
since it is
obviously met by all proximity measures under consideration).

\statemento{Independence of other components.}{Let $A$ and $B$ be two
different components of a multigraph. Then any addition,
removal, or reweighting of edges $($arcs$)$ within $B$ does not alter
the values of proximity for the vertices that belong to $A$.  }

\PRP{5}{
Relative forest accessibility for multigraphs has the following
properties$:$
symmetry$,$
nonnegativity$,$
diagonal maximality$,$
the triangle inequality for proximities$,$
the disconnection condition$,$
the transit property$,$
monotonicity$,$
doubly stochastic normalization$,$ and
macrovertex independence.
}

Thereby, relative forest accessibility for multigraphs possesses all
normative properties of Sec.~2
without any restrictions on the
weights of edges, and it features macrovertex independence and doubly
stochastic normalization. Certainly, this does not raise relative forest
accessibility over other proximity measures. Rather, this index perfectly
corresponds to one possible concept of proximity specified by the properties
listed in Proposition~5.

\arcsect{7}{COMPONENTS OF RELATIVE FOREST ACCESSIBILITY}

In this section, the relative forest accessibility for multigraphs is
decomposed into components that correspond to the sets of forests with a
varying number of trees. Next, we consider the notions of proximity that
correspond to each component. Let $v$ be the number of connected components in
$G$; by $V_i$ we denote the set of vertices of the component of $G$ that
contains vertex $i$ ($i=\1n$).

\statement{THEOREM~2~(a parametric version of the matrix-forest theorem for
multigraphs).}
{For any weighted multigraph $G$ and any $\tau\ge0$, let
$Q(\tau)=(q\_{ij}(\tau))$ be the matrix $(I+{\tau} L)^{-1}.$ Then $Q(\tau)$
exists and
\begin{equation}
q\_{ij}(\tau)=\sum^{n-v}_{k=0}\tau^k{\e}(\FF^{ij}_k)
\Big/ \sum^{n-v}_{k=0} \tau^k{\e}(\FF_k),\quad i,j=\1n,
\label{f1}
\end{equation}
where $\FF_k$ is the set of spanning rooted forests in $G$ that consist of $k$
edges, and $\FF^{ij}_k\subseteq\FF_k$ is its subset comprising those forests
in which $j$ belongs to a tree rooted at $i$.
}

By Proposition~5, the matrix of relative forest accessibilities is doubly
stochastic, whence $\suml^n_{j=1}q\_{ij}(\tau)=1$, $i=\1n$, $\:\tau\ge0$. The
following proposition states a stronger fact, namely, the stochastic property
is true for the coefficients at every exponent of $\tau$ in~(\ref{f1}).

\PRP{6}{
For any $i=\1n$ and $k=0\cdc n-v$, we have
\begin{equation}
\sum^n_{j=1}\e(\FF^{ij}_k)=\e(\FF\_k).
\label{rowsum}
\end{equation}
}

The matrices $Q(\tau)$, $\:\tau>0$, make up a parametric family of relative
forest accessibility indices which obviously have the same basic properties
as $Q=Q(1)$. By (\ref{f1}), $Q(\tau)$ can be represented as
\begin{equation}
Q(\tau)={1 \over s(\tau)} \left({\tau}^0 Q_0+{\tau}^1 Q_1+\ldots+{\tau}^{n-v}
Q_{n-v}\right), \label{expa}
\end{equation}
where $s(\tau)=\suml^{n-v}_{k=0}{\tau}^k\e(\FF_k)$, $Q_k=(q\_{k,ij})$, and
$q\_{k,ij}=\e(\FF^{ij}_k),\;k=0\cdc n-v,\>i,j=\1n$.

Every matrix $Q_k,\;$ $k=0\cdc n-v,$ reflects a specific vertex proximity. Let us consider them in some
detail. First, $Q_0=I$, i.e., the ``proximity'' specified by $Q_0$ is simply identity. Further, the entry
$q\_{1,ij},$ $j\ne i,$ of $Q_1$ is equal to the total weight of the edges in $G$ that are incident to $i$ and
$j$. Generally, the entry $q\_{k,ij}$ of $Q_k$ is distinct from zero if and only if $G$ contains some paths
of length $k$ or shorter between $i$ and $j$. The corresponding notion of proximity ignores all paths of
length $k+1$ or longer. Whenever $k\ge|V_i|_{\max}-1$ (where $|V_i|_{\max}$ is the maximum number of vertices
among the components of $G$), the proximity corresponding to $Q_k$ takes into account all paths in $G$.

Recall that $V_i$ is the set of vertices in the component of $G$ that contains
$i$. To examine the proximity corresponding to $Q_{n-v}$, we
introduce the matrix $\oj (G)=\oj=(\oj_{ij})$:
\[
{\oj}_{ij}
=\left\{
\begin{array}{ll}
{{\textstyle{1}}\over\textstyle{|V_i|}}, & \mbox{if $j\in V_i$,} \\
0,                                       & \mbox{otherwise}
\end{array}
\right.
\]
and prove the following lemma.

\LE{1}{
\begin{equation}
Q_{n-v}=\e(\FF_{n-v})\oj.
\label{stupid}
\end{equation}
}

As mentioned above, the ``proximity'' that corresponds to $Q\_0$ is identity. By Lemma~1, the matrix
$Q_{n-v}$ represents an opposite concept of proximity: all vertices that belong to the same component of $G$
are equally ``close'' to each other, and the value of their proximity is inversely proportional to the number
of vertices in the component. Thus, the proximity to vertex $i$ is {\em uniformly distributed\/} over the
component of $G$ that contains $i$. If $G$ is connected, then $\oj=(1/n)J$, where $J$ is the $n\times
n$-matrix having all entries one, and so all entries of $Q_{n-v}$ are $\e(\FF_{n-v})/n$. For all matrices
$Q_k$, $k=0\cdc n-v,$ the proximity of two vertices from different components of $G$ is zero.

\CRL{1}{
$\liml_{\tau\to\infty} Q (\tau)=\oj.$
}

Corollary~1 follows directly from Theorem~2 and Lemma~1.

\RM{1}
The matrix $Q_{n-v-1}$ is of special interest. Its entry
$q\_{n-v-1,ij}$ is the total weight of those spanning rooted forests in $G$
that

(1) have two trees in one component of $G$ and one tree in each of the others,
and

(2) have $i$ and $j$ in the same tree rooted at $i$.

Among the matrices $Q_k,$ $k=0\cdc n-v,$ the matrix $Q_{n-v-1}$ is the most
similar (in the properties) to the matrices $Q(\tau)$ of relative forest
accessibility. Indeed, by (\ref{expa})--(\ref{stupid}), the comparison of
two entries of $Q(\tau)$ at a large $\tau$ is determined by the comparison of
the corresponding entries of $Q_{n-v-1}$. Only when the latter entries are
equal do the corresponding entries of $Q_k$ at $k<n-v-1$ matter.
Dealing with examles convinces us that the situations where
two entries of $Q_{n-v-1}$ are equal, whereas the corresponding entries
of $Q_k,$ $k<n-v-1,$ vary, are not frequent, and indeed  it is not easy
to intuitively discriminate between the compared proximities in these
cases. Still, an important exception exists. As mentioned above,
$k\ge|V_i|_{\max}-1$ is necessary and sufficient for $Q_k$ to take into
account all paths in $G$. If all components of $G$, except one, are
separate vertices or $G$ is connected, then $|V_i|_{\max}-1=n-v$. In this
case, if a pair of vertices in the nontrivial component is connected only by
paths of length $n-v$ (a chain graph), then the corresponding entry of
$Q_{n-v-1}$ is zero, and so $Q_{n-v-1}$ violates the disconnection condition.
Note that some weighted sums of $Q_{n-v-1}$ and $Q_{n-v}$ are free of this
flaw.  Such linear combinations are studied in the following section.
Moreover, we show that $Q_{n-v-1}$ is closely connected with the matrix
$\Lp$, the Moore--Penrose generalized inverse of $L$. More precisely,
$\Lp$ is the sum of $Q_{n-v-1}$ and $Q_{n-v}$ with definite
coefficients.

\arcsect{8}{ACCESSIBILITY VIA DENSE FORESTS CONNECTED WITH THE
GENERALIZED INVERSION OF THE LAPLACIAN MATRIX}

This section is devoted to weighted sums of matrices $Q_{n-v-1}$ and
$Q_{n-v}=\e(\FF_{n-v})\oj$. A number of papers \cite{KleRan,Kle97,Altm,SteZel}
use, either explicitly or implicitly, proximity matrices whose generalization
to multicomponent graphs can be represented as ${(L+\alpha\oj)}^{-1}$, where
$\alpha>0$. The aims of this section are as follows:

(1) to provide a topological interpretation of such a proximity in the case of
arbitrary multigraphs (it is based on the matrices $Q_{n-v-1}$ and $Q_{n-v}$);

(2) to establish its relation with the matrix $\Lp$, the Moore--Penrose
generalized inverse of $L$, and

(3) to ascertain its properties.

We will show that ${(L+\alpha\oj)}^{-1}$
with a sufficiently small $\alpha$ is a weighted sum of $Q_{n-v-1}$ and
$Q_{n-v}$ with positive coefficients and satisfies a number of conditions
of Sec.~2.

To solve the foregoing problems, we will need the matrix
\begin{equation}
\widetilde{Q}={(L+\oj)}^{-1}-\oj,
\end{equation}
which has many remarkable properties. Four representations for $\widetilde{Q}$
are stated below (Proposition~7--9 and Theorem~3).

\PRP{7}{
For any
$\alpha\ne0,$ the matrix $(L+\alpha\oj)$ is invertible, and
$\widetilde{Q}=(L+\alpha\oj)^{-1}-\alpha^{-1}\oj.$}

By Proposition~7, the difference between $\widetilde {Q}$ and
$(L+\alpha\oj)^{-1}$ is represented by a matrix whose entries are constant
within each component of $G$. In \cite{KleRan,Kle97,Altm,SteZel},
matrices of the form of $(L+\alpha\oj)^{-1}$ are mainly used for
transformations such as (\ref{metr}), where, if one pays no regard for
intercomponent entries, they can be equivalently replaced by
$\widetilde{Q}$.

Recall that for any rectangular complex matrix $A$, the Moore--Penrose
generalized inverse of $A$ is the unique matrix $\Ap$ such that

(1) $A\Ap$ and ${\Ap}A$ are Hermitian matrices,

(2) $A{\Ap}A=A$, and

(3) $\Ap{A}\Ap=\Ap$.

\PRP{8}{
For any weighted multigraph $G$, the matrix $\widetilde{Q}$ is the
Moore--Penrose generalized inverse of $L=L(G)$, that is,
$\widetilde{Q}=\Lp$.  }

Since $L$ is a square matrix, and $A\Ap={\Ap}A$ (which follows from the
proof of Proposition~8), the matrix $\widetilde{Q}$ is the {\it group
inverse\/} of $L$ (cf.~\cite{KirNeuSha}). Geometric interpretations for
$\Lp$ are given in~\cite{Fiedler95}.

It turns out that $\Lp$ can be obtained by a passage to the limit from the
parametric matrix $Q(\tau)$ of relative forest accessibilities
(cf.~Corollary~1).

\PRP{9}{
$\Lp=\liml_{\tau\to\infty}\tau(Q(\tau)-\oj).$
}

Proposition~9 and Theorem~2 enable one to obtain a topological interpretation
for $\Lp=(\lpij)$.

\statement{THEOREM~3~(a topological interpretation for the matrix $\Lp,$
the Moore--Penrose generalized inverse of $L$):}
{
\begin{equation}
\lpij=\cases{
{{\textstyle\e(\FF^{ij}_{n-v-1})-{1\over|V_i|}\cdot\e(\FF_{n-v-1})}
\over{\textstyle\e(\FF_{n-v})}},             &if $\>j\in V_i$,\cr
0,                                           &otherwise.\cr}
\label{tilQF}
\end{equation}
}

Here, the numerator is the result of centralization: the $ij$-entry minus the
$i$th-row mean of $Q_{n-v-1}$ (see~(\ref{rowsum})). By Theorem~3,
the definition of $\oj$, and Lemma~1, one has
\begin{eqnarray}
\Lp
&=&{\e(\FF_{n-v-1})\over\e(\FF_{n-v})}
\left({1\over\e(\FF_{n-v-1})}Q_{n-v-1}
-     {1\over\e(\FF_{n-v  })}Q_{n-v}
\right)
\nonumber\\
&=&{1\over\e(\FF_{n-v})}
\left(Q_{n-v-1}-\e(\FF_{n-v-1})\oj\right).
\label{Qkruto}
\end{eqnarray}

Another representation of $\Lp$ for connected weighted graphs was obtained
in~\cite{KirNeuSha}.

Can $\Lp$ be considered as a matrix of vertex proximities? By (\ref{tilQF}), this ``proximity'' equals zero
for vertices from different component of $G$, and so does the sum of ``proximities'' of each vertex with the
vertices of the same component. The latter does not match an intuitive idea of proximity.  First,
nonnegativity is violated; second, the ``proximity'' of poorly connected vertices from the same component
turns out to be less than that for any vertices from different components.

Now, let us return to the matrices $(L+\alpha\oj)^{-1}.$
Propositions~7--9 and Eq.~(\ref{Qkruto}) imply the following identities:
\begin{eqnarray}
(L+\alpha\oj)^{-1}&=&\Lp                       +\alpha^{-1}\oj\\
                  &=&\liml_{\tau\to\infty}\tau(Q(\tau)-\oj)
                                               +\alpha^{-1}\oj
\label{QalphaLim}\\
                  &=&{1\over\e(\FF_{n-v})}
                     \left(Q_{n-v-1}+\left(\alpha^{-1}-
                     {\e(\FF_{n-v-1})\over\e(\FF_{n-v})}
                     \right)Q_{n-v}\right)
\label{QalphaQQ}\\
                  &=&{1\over\e(\FF_{n-v})}
                     Q_{n-v-1}+\left(\alpha^{-1}-
                     {\e(\FF_{n-v-1})\over\e(\FF_{n-v})}
                     \right)\oj.
\label{Qalpha}
\end{eqnarray}

Thus, whenever $0<\alpha<{\e(\FF_{n-v})/\e(\FF_{n-v-1})}$, the matrix
$(L+\alpha\oj)^{-1}$ is the sum of $Q_{n-v-1}$ and $Q_{n-v}$ with 
positive coefficients. Let a {\em dense forest\/} be a
spanning rooted forest in $G$ with $n-v$ or $n-v-1$ edges. Then the proximity
measure (\ref{QalphaQQ}) with $0<\alpha<{\e(\FF_{n-v})/\e(\FF_{n-v-1})}$
can be referred to as {\it accessibility via dense forests}.

\PRP{10}{
The accessibility via dense forests in the case of multigraphs has the
following properties$:$
symmetry$,$
nonnegativity$,$
diagonal maximality$,$
the triangle inequality for proximities$,$
the disconnection condition$,$ and
the transit property$.$
It does not satisfy monotonicity.
}

It is interesting to examine the nature of the violation of monotonicity.
It follows from (\ref{QalphaLim}) that whenever $k$ and $t$ belong to the same
component of the original multigraph, monotonicity is valid in a
nonstrict form, i.e., all strict inequalities are replaced by nonstrict ones,
which can be regarded as acceptable. Rough violations of monotonicity (namely,
$\D p\_{kt}<\D p\_{ij}$ and $\D p\_{kt}<0$) only occur when $k$ and $t$
originally belong to different components of $G$. This suggests an idea of
searching for a better modification of accessibility via dense forests.
The scrutiny of this question, as well as the examination of the metric
corresponding (in the sense of \cite{weAiT98}) to this proximity measure (see
\cite{KleRan,Kle97,Altm,KirNeuSha}), is beyond the scope of this paper.

\arcsect{9}{ON SOME PECULIARITIES OF THE PROXIMITY MEASURES}

A specific feature of path and route accessibilities is the necessity
of imposing rather strong restrictions on the weights of edges (arcs)
to guarantee the properties of Sec.~2 
convergence (in the case of route accessibility). These restrictions imply a fast decrease of proximity with
movement away from a vertex along an edge chain. A characteristic feature of connection reliability is the
effect of saturation. If, for example, two vertices are connected by an edge, the weight of which is close to
1, then the addition of other paths between them leaves the value of proximity almost the same.  In addition,
all diagonal entries are ones, i.e., they do not characterize self-relations of any kind. Accessibility via
dense forests violates monotonicity when two components of a graph get a connection; it only satisfies the
nonstrict version of monotonicity, when a graph is changed within components. Unlike relative forest
accessibility, here the triangle inequality is also satisfied in a nonstrict form, provided that $i,j$, and
$k$ are distinct. On the other hand, the metric derived from this proximity measure by (\ref{metr}) coincides
with the classical graph metric in the case of trees~\cite{KleRan}.  For a further study of this metric,
see~\cite{KirNeuSha}. The relative forest accessibility differs from the other proximity measures by the very
fact of its relativeness. A manifestation of this is the stochastic normalization property of the matrices
$Q$ and $Q(\tau)$ for digraphs and doubly stochastic normalization in the case of undirected graphs. As a
corollary, the addition of new edges (arcs) in a graph does not increase all proximities; some of them will
necessarily decrease. The corresponding ``absolute'' proximity measure can be obtained by considering the
adjugate of the matrix $(I+\tau L)$ instead of $Q(\tau)=(I+\tau L)^{-1}$. In addition, relative forest
accessibility features macrovertex independence, which is not always desirable. To illustrate these and some
other peculiarities of the proximity measures under study, we shall consider a few simple examples.

\begin{figure}[t]
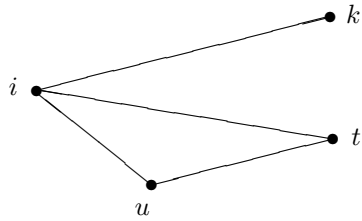

\bigskip
\input p2a.pic
\vspace{-2em}
\caption{Example 1.}
\label{P_sig}
\end{figure}

For the graph in Fig.~\ref{P_sig}, path accessibility, connection reliability,
and route accessibility give $p\_{ik}<p\_{it}$.
Seemingly, it would otherwise be unnatural, since $i$ and $t$ are connected not only by an edge (as $i$ and
$k$ are), but also by a path of length 2 $(iut)$. Nevertheless, the relative forest accessibility gives
$p\_{ik}=p\_{it}=p\_{iu}$ (this follows from macrovertex independence: $\{k,t,u\}$ is a macrovertex). The
same result is provided by the accessibility via dense forests. Macrovertex independence is appropriate when
any connections within a macrovertex can be regarded as its ``domestic affairs.'' For example, if each
professor gives his/her lectures to all students (then the students form a macrovertex), and the students
write them down verbatim, then no reading or rewriting of the notes of each other can help them learn
anything more (i.e., to approach the knowledge of the professors).

\begin{figure} [t]
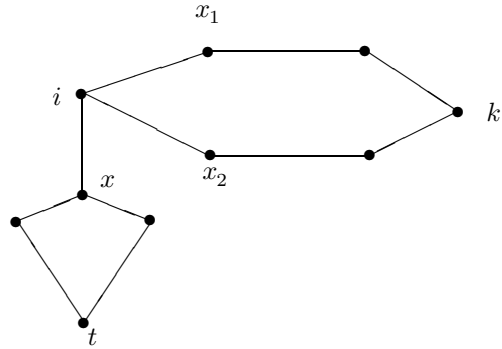

\input p3b.pic
\vspace{-2em}
\caption{Example 2.}
\label{P_put}
\end{figure}

The following example illustrates some peculiarities of the path and route
accessibilities. In Fig.~\ref{P_put}, $i$ is connected with $k$ by
two paths, as well as with $t$, and the weights of these paths are equal
(provided that the weights of all edges are equal). Hence, the path
accessibilities $p\_{ik}$ and $p\_{it}$ are also equal. But the paths that
connect $i$ to $t$ have a common edge.
Therefore, connection reliability gives $p\_{ik}>p\_{it}$. The same result
holds for relative forest accessibility and accessibility via dense forests.
In contrast, route accessibility provides $p\_{ik}<p\_{it}$. This is because
there exist two paths of length two from $x$ to $t$ and only one path of
length two from $x\_1$ (or from $x\_2$) to $k$. As a result, there are eight
routes of length seven from $i$ to $t$ and only four routes of length seven
from $i$ to $k$.

\begin{figure}[t]
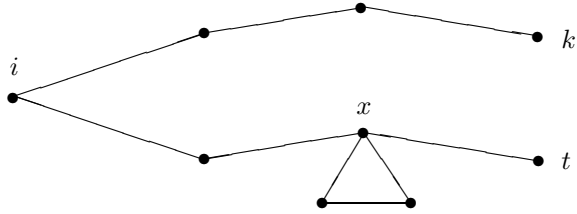

\input p4.pic
\vspace{-2em}
\caption{Example 3.}
\label{P_cyc}
\end{figure}

Furthermore, the proximity measures at hand behave differently as
applied to cycles. The cycle in Fig.~\ref{P_cyc} has no influence on
the values of path accessibility and connection reliability between $i$
and $t$, i.e., $p\_{it}=p\_{ik}$ (if all edge weights are equal).
Using route accessibility, we have $p\_{it}>p\_{ik}$. At the same time,
relative forest accessibility provides $p\_{it}<p\_{ik}$, as the
approach of $i$ and $t$ to the vertices of the cycle (owing to its appearance)
moves them away (in the relative account) from each other. The same holds for
the accessibility via dense forests.
\begin{table}[t]
\begin{center}
\begin{tabular}{l}
\begin{tabular}{lccccccc}
\hline
Property                  &Paths&Reliability& Routes &\begin{tabular}{l}
                                                      Forests\\
                                                      (undirected)
                                                      \end{tabular}
                                                                  &\begin{tabular}{l}
                                                                   Dense forests\\
                                                                   (undirected)
                                                                   \end{tabular} \\
\hline
Symmetry                  &$+$  & $+$       &$+$     &$+$         &$+$           \\
Nonnegativity             &$+$  & $+$       &$+$     &$+$         &$+$           \\
Reversal property         &$+$  & $+$       &$+$     & $\times$   & $\times$     \\
Diagonal maximality       &\ps  & \ps       &$+$     &$+$         &$+$           \\
\begin{tabular}{l}%
$\!\!\!\!$Triangle inequality\\
for proximities
\end{tabular}             &\ps  & \ps       &\pb     &$+$         &$+$           \\
Disconnection condition   &$+$  & $+$       &$+$     &$+$         &$+$           \\
Transit property          &\ps  & \ps       &$+$     &$+$         &$+$           \\
Monotonicity (1)          &\ps  & \ps       &$+$     &$+$         & $-$          \\
Monotonicity (2)          &$+$  & \ps       &$+$     &$+$         & $-$          \\
Monotonicity (3)          &$+$  & $+$       &$-$     &$+$         & $-$          \\
\hline
\end{tabular}\\
\noalign{\vskip 0pt}
$+^*$\ is valid under some restriction and/or in the nonstrict form.\\
$+^{**}$\ was proved under an additional constraint.\\
$\times$ inapplicable, as only undirected graphs were considered.
\end{tabular}\\
\caption{Some properties of proximity measures for graph vertices.}
\label{Table1}
\end{center}
\end{table}

Note finally that for path accessibility, connection reliability, and the
measures representable by weighted sums of the matrices $Q_1\cdc Q_{n-v}$
with fixed weights, the values of proximity linearly depend on the weights of
edges (arcs), whereas for the other measures at hand, this is not the case.

Thus, the proximity measures under discussion have significantly different properties. At the same time,
``almost all'' of them possess ``almost all'' of the ``basic'' properties formulated in Sec.~2

\arcsect{10}{CONCLUSION}

In this paper, we have dealt with several proximity measures for the
vertices of directed and undirected multigraphs and considered their
properties. These properties and the informal discussion of the
previous section can help one choose adequate proximity measures when
exact mathematical models are lacking.

A common feature of the indices considered in this paper is the measurement of the proximity (accessibility,
connectivity) of two vertices by the total weight of certain substructures that ``connect'' these vertices.
As such substructures, we examined paths (in particular, taking into account their overlaps), routes,
spanning rooted forests, and ``dense'' spanning rooted forests. The weight of a substructure was defined as
the product of the weights of the constituent edges (arcs). Within this approach, a proportional modification
of all edge weights is needed in some cases, as well as assigning the same weight to all edges (arcs) of
unweighted graphs. In conclusion, let us indicate some proximity measures that do not enter into the scope of
the present paper.  These are the indices dual (in the sense of~\cite{weAiT98}) to the classical distance for
connected graphs and to some nonclassical graph metrics \cite{Kle97}, maximum flow (minimum cut) between
vertices \cite{Altm}, and a number of measures related to random walks in graphs
(see~\cite{DoyleSnell,Altm,Kle97}).

\arcpril 

\PPR{1}
{\it Symmetry}, {\it nonnegativity}, {\it
reversal property}, and {\it disconnection condition\/} immediately follow
from the definition of path accessibility. To prove the remaining properties,
let us find $\e\_0$ guaranteeing that whenever the weights of all edges (arcs)
are less than $\e\_0$, $p\_{ij}<1$ holds for all $i$ and $j\ne i$. Let $m$ be
the greatest possible number of edges (arcs) incident to the same pair of
vertices. Note that when a multigraph $G$ is complete (i.e., exactly $m$
edges are incident to each pair of vertices), and the weights of all edges are
$\e$, then at $j\ne i,\;$ $p\_{ij}=\suml_{k=1}^{n-1}A_{n-2}^{k-1}(\e m)^k$,
where $A_{n-2}^{k-1}$ is the number of permutations of $n-2$ things taken
$k-1$ at a time. Now, we equate this expression to unity and assign to
$\e\_0$ the positive root of the equation obtained.

Henceforth, we will assume that $\e_{ij}^p<\e\_0$ for all edge weights
$\e_{ij}^p$.  As $p\_{ij}$ is maximal in a complete multigraph, this
will guarantee \begin{equation} p\_{ij}<1,\quad i,j=\1n,\; i\ne j,
\label{modest}
\end{equation}
for all weighted multigraphs on $n$ vertices with the number of multiple edges
not greater than $m$. The same constraint can be obtained for the weights of
arcs in multidigraphs.

{\it Diagonal maximality\/} follows from the inequalities $p\_{ii}\ge1$ and
(\ref{modest}).

Prove the {\it triangle inequality for proximities}. At $i=j$ or
$i=k$, the inequality reduces to equality. Suppose that $i\ne j$ and $i\ne k.$
Note that whenever all paths from $j$ to $k$ pass through $i$, $p\_{jk}
=p\_{ji} p\_{ik}$ holds; otherwise $p\_{jk}\ge p\_{ji}p\_{ik}$. Let
$C$ be the total weight of simple cycles from $i$ to $i$; then $p\_{ii}=1+C$.
Using (\ref{modest}), one obtains the triangle inequality for proximities:
\[
p\_{ij}+p\_{ik}-p\_{jk}-p\_{ii}\le
p\_{ji}+p\_{ik}-p\_{ji}p\_{ik}-1-C=(p\_{ji}-1)(1-p\_{ik})-C<0.
\]

To prove the {\it transit property}, note that $p\_{it}=p\_{ik}p\_{kt}$,
and using (\ref{modest}), we have $p\_{it}<p\_{ik}.$

Now prove {\it monotonicity}. Item~1. Suppose that $\D\e\_{kt}$ is the increment of the weight of an existing
edge or the weight of a new edge between $k$ and $t$. Then $\D p\_{kt}= \D\e\_{kt}>0$. Let us show that
whenever all edge weights are smaller than $\e\_0$ and $\{i,j\}\ne\{k,t\}$, $\D p\_{ij}<\D\e\_{kt}$ holds. If
$i=k$, then $\D p\_{ij}=\D p\_{kj}\le\D\e\_{kt}p\_{tj}$, and the required inequality follows from
(\ref{modest}). The cases $i=t$, $j=k$, and $j=t$ are similar. It remains to consider the case
$\{i,j\}\cap\{k,t\}=\varnothing$, in which $n\ge4$. Obviously, $\D p\_{ij}=\D\e\_{kt}w$, where $w$ is the
total $\overline{(k,t)}$-weight of the paths from $i$ to $j$ that contain the new (reweighted) edge $(kt)$,
and the ``$\overline{(k,t)}$-weight'' of a path is the product of the weights of all its edges, except for
the edge $(kt)$. Prove that $w<1$. Obviously, $w$ is maximal in a complete multigraph, where, as is easy to
check, $w=2\suml_{k=2}^{n-2}(k-1)A_{n-4}^{k-2} (\e\_0 m)^k$. Let us show that in this case, $w$ is less than
the value $p=\suml_{k=1}^{n-1}A_{n-2}^{k-1}(\e\_0 m)^k$ of the proximity for two distinct vertices in a
complete multigraph, which equals 1 by the definition of $\e\_0$. Juxtapose the coefficients at the same
exponents of $(\e\_0 m)$ in the expressions for $w$ and $p$. It is easy to verify that the inequality
$2(k-1)A_{n-4}^{k-2}\ge A_{n-2}^{k-1}$ has a unique solution: $n=4,$ $k=2$. Thereby, the statement is proved
in the case of $n>4$. Finally, for $n=4$ we have $p=\e\_0 m+2(\e\_0 m)^2+2(\e\_0 m)^3$ and $w=2(\e\_0 m)^2$;
therefore, $w<p$ as well. A similar proof applies to multidigraphs.

Item~2. The statement follows from $\D p\_{ik}=0$ and $\D p\_{it}>0.$

Item~3. We have $\D p\_{i\_1 i\_2}=0$, as the edge (arc) $(kt)$ does
not belong to any path from $i\_1$ to~$i\_2$.
\epr

\PPR{2}
{\it Symmetry}, {\it nonnegativity}, {\it reversal property}, and {\it
disconnection condition\/} follow easily from the definition of connection
reliability. {\it Diagonal maximality\/} in a nonstrict version follows from
the facts that $p\_{ii}=1$ and $p\_{ij}\le1,\:$ $i,j=\1n$. If all edge/arc
weights are less than~1, then, obviously, $p\_{ij}<1$ at $j\ne i$;
therefore $p\_{ij}<p\_{ii}$.

The proof of the {\it triangle inequality for proximities\/} mimics the
corresponding proof for path accessibility.

{\it Transit property\/} (in the form specified in Proposition~2) follows
from the equality $p\_{it}=p\_{ik} p\_{kt}$, which is valid under the
hypothesis of this property.

Prove item~1 of {\it monotonicity\/} for multidigraphs. This proof will
also be applicable to multigraphs. Let a {\em state\/} of a
multidigraph, all of whose arcs are assigned some intactness
probabilities, be any of its spanning subgraphs.  The arcs of the
subgraph are interpreted as the only intact arcs of the original
multidigraph.  By the assumption of independence of failures,
the probability of a state is the product of the intactness
probabilities of the arcs entering into the state and the failure
probabilities of the lacking arcs. Let a new arc from $k$ to $t$ be
added. Note that $\D p\_{ij}$ is the total probability of those states
in which

(1) the new arc $(kt)$ is present,

(2) there is a path from $i$ to $j$, and

(3) the removal of the arc $(kt)$ leaves no path from $i$ to $j$.

Note that in all these states, the removal of $(kt)$ does not leave any path
from $k$ to $t$ either (otherwise the removal of this arc would not have
broken a path from $i$ to $j$). Therefore, the specified total probability is
a summand of $\D p\_{kt}$, and hence, $\D p\_{kt}\ge\D p\_{ij}$. Whenever all
arc weights are strictly less than~1, there is at least one state whose
nonzero probability is a summand of $\D p\_{kt}$, but does not enter into $\D
p\_{ij}$: in this state the new arc $(kt)$ is solely intact, and the desired
inequality is strict. All these conclusions are preserved when the weight of
an arc $(kt)$ {\em increases}. This is because the connection reliability is
{\em affinely related\/} with each arc weight.

Item~2 of {\it monotonicity\/} is true, since $\D p\_{ik}=0$, $\D p\_{it}\ge0$,
and $\D p\_{it}>0$ when all arc/edge weights are strictly less than one. Item~3
is valid, as the edge (arc) $(kt)$ does not belong to any path from $i\_1$
to~$i\_2$.
\epr

\PPR{3}
{\it Symmetry}, {\it nonnegativity}, {\it
reversal property}, and {\it disconnection condition\/} follow from the
definition of route accessibility.

Prove {\em diagonal maximality\/} for multidigraphs.
In talking about route accessibility, we always consider a family of
graphs with a specified greatest possible number of multiple edges (arcs) $m$
and with edge/arc weights smaller than $\e\_{\max}=(m(n-1))^{-1}$. Suppose
that $\G$ is a weighted multidigraph that belongs to such a family; $i$ and
$j\ne i$ are arbitrary vertices of $\G$; $\e<\e\_{\max}$ is the maximum among
the arc weights in $G$. Consider the multidigraph $\G'$ constructed by
removing all arcs directed to $i$ from the complete multidigraph with the
multiplicity of all arcs $m$ and the weight of all arcs $\e$. Obviously, for
$\G',$ $p'_{ii}=1$, and $p'_{ij}=p'_{ik}$ for any $k\ne i$. Particularizing
the equality $P'(I-E')=I$ for the $ij$-entry of $P'(I-E'),$ we derive
$p'_{ij}=\frac{\textstyle\e m}{\textstyle 1-(n-2)\e m}$, and consequently,
$p'_{ii}>p'_{ij}$, since $\e<(m(n-1))^{-1}$. If some arcs are removed from
$\G'$ or some weights of arcs are reduced (let the resulting graph be $\G''$),
then $p\_{ii}$ does not change, whereas $p\_{ij}$ can only decrease. Now, let
an arc from $k\ne i$ to $i$ be added to $\G''$. By virtue of Proposition~4
(the proof of which is given below), in this case $\D p\_{ii}-\D p\_{ij}=
hp''_{ik}(p''_{ii}-p''_{ij})>0$, and thus, $p\_{ii}>p\_{ij}$ remains
true.  Similarly, $p\_{ii}>p\_{ij}$ is preserved at the consecutive
addition of other arcs directed to $i$. Hence, $p\_{ii}>p\_{ij}$ is
also valid for $\G$, and the diagonal maximality is proved. The
fulfillment of this property for any multigraph $G$ is ensured by its
validity for the symmetric multidigraph $\G$ with the same matrix $E$.

Now we prove {\em triangle inequality for proximities\/} in the case
where the weights of all edges (arcs) do not exceed $(mn)^{-1}.$ First,
consider the digraph $\G'$ that differs from the complete digraph by the lack
of all arcs directed to $i$. At $j=i$ or $k=i$, the triangle inequality for
proximities reduces to equality, so assume that $j\ne i$ and $k\ne i$. Let
each arc of $\G'$ have weight $\e=1/n$. Using the equality $(I-E')P'=I$ for
the entries $ij,$ $ik$, and $ii$ of $(I-E')P'$, one obtains
\begin{eqnarray*}
p'_{ij}&=&p'_{ik}=\frac{\e}{1-(n-2)\e}={1\over 2},\\
p'_{ii}&=&1,
\end{eqnarray*}
hence, $p'_{ii}-p'_{ij}-p'_{ik}+p'_{jk}>p'_{ii}-p'_{ij}-p'_{ik}=0.$ We shall
prove now that no change of $\G'$ can decrease $p\_{ii}-p\_{ij}-p\_{ik}.$
Indeed, if some arcs are removed from $\G'$ and/or the weights of some
arcs are reduced, $p\_{ii}$ does not change, whereas $p\_{ij}$ and
$p\_{ik}$ can only decrease; therefore,
$p\_{ii}-p\_{ij}-p\_{ik}\ge0$ is preserved. Furthermore, if for some
digraph $\G$ this inequality is valid, then the addition of any arc
$ti$ to $\G$ cannot violate it, since, by Proposition~4, \[ \D
p\_{ii}-\D p\_{ij}-\D
p\_{ik}=h(t)\,p\_{it}(p\_{ii}-p\_{ij}-p\_{ik})\ge0.  \] Thus, the
triangle inequality for proximities is valid for any digraph. The
fulfillment of this property for multidigraphs is proved by replacing the set
of arcs between a pair of vertices with a single arc with the total
weight, which reduces the problem to digraphs. The fulfillment of the property
for any multigraph is ensured by its validity for the symmetric multidigraph
with the same matrix $E$.

{\it Transit property\/} for multidigraphs will be proved by
contradiction. Let $\G$ be the multidigraph with the minimum number of
arcs among the multidigraphs that violate the transit property. Then
$\G$ has a path from $i$ to $k$, $t\ne k$, and any path from $i$ to $t$
contains $k$, but $p\_{ik}\le p\_{it}$. From the diagonal maximality,
$k\ne i$. Let $(ij)$ be the first arc of an arbitrary path from $i$ to $k$,
and let $\G'$ be the multidigraph obtained by removing the arcs $(ij)$ from
$\G$. Then, after adding the arc $(ij)$ to $\G'$, one has $\D p\_{it}
\ge\D p\_{ik}$. Indeed, if $\G'$ has no path from $i$ to $k$, then
$p'_{ik}=p'_{it}=0$ in $\G'$, and $\D p\_{it}<\D p\_{ik}$ would have been in
contradiction with $p\_{ik}\le p\_{it}$ in $\G$. If, otherwise, $\G'$ contains
a path from $i$ to $k$ and $\D p\_{it}<\D p\_{ik}$, then $\G'$ violates the
transit property, which contradicts the minimality of $\G$. Further, by
Proposition~4, $\D p\_{it}-\D p\_{ik}=hp'_{ii}(p'_{jt}-p'_{jk})$, where $h>0$,
and $\D p\_{it}\ge\D p\_{ik}$ implies $p'_{jt}\ge p'_{jk}$. By the
construction, $\G'$ has a path from $j$ to $k$, and any path from
$j$ to $t$ contains $k$. Hence, $\G'$ breaks the transit property, which
contradicts the minimality of $\G$. {\it Transit property\/} for any
multigraph is proved by turning to the multidigraph with the same matrix
$E$.

To prove item~1 of {\it monotonicity\/} in the case of multidigraphs, note
that, by virtue of Proposition~4, $\D p\_{kt}=hp\_{kk} p\_{tt}$ and
$\D p\_{ij}=hp\_{ik}p\_{tj}.$ Now, the required statement follows from the
diagonal maximality and can be extended to multigraphs by a standard trick.
Similarly, item~2 of monotonicity follows from the formula $\D p\_{it}-
\D p\_{ik}=h(p\_{ik}p\_{tt}-p\_{ik}p\_{tk})$ and diagonal maximality.
Item~3 is not true, since under the hypothesis of monotonicity, some routes
from $i\_1$ to $i\_2$ that contain the edge (arc) $(kt)$ can appear or
increase their weight.
\epr

\PPR{4}
Let $\D (I-E)=(I-E')-(I-E).$ Note that
$\D(I-E)=XY$, where $X=(x\_{i1}),\:$ $i=\1n,$ is the column vector with entries
$x\_{k1}=-\D\e\_{kt}$ and $x\_{i1}=0$ for all $i\ne k$; $Y=(y\_{1j}),\:$
$j=\1n$, is the row vector with entries $y\_{1t}=1$ and $y\_{1j}=0$ for all
$j\ne t.$ According to \cite[Sec.~0.7.4]{HoJo},
\[
P'=P-\frac{1}{1+YPX}PXYP.
\]
It is straightforward to verify that
$(-\frac{1}{1+YPX})=-h/\D\e\_{kt}$ and $\:PXYP=-\D\e\_{kt}R$, and thereby the
proposition is proved.
\epr

\PPR{5}
Let us prove item~1 of {\em monotonicity\/} (all the other statements are
proved in \cite{weAiT97}). By item~1 of Proposition~7 from \cite{weAiT97},
$\D p\_{kt}=h(p\_{kk}-p\_{kt})(p\_{tt}-p\_{tk})$ and
$\D p\_{ij}=h(p\_{ik}-p\_{it})(p\_{jt}-p\_{jk})$, where $h>0$. Diagonal
maximality implies $\D p\_{kt}>0$. If $\D p\_{ij}>0$, then
$(p\_{ik}-p\_{it})(p\_{jt}-p\_{jk})>0$. For definiteness, we assume
that $p\_{ik}-p\_{it}>0$ and $p\_{jt}-p\_{jk}>0$ (the complementary
case is treated similarly). Then, by item~2 of Proposition~6 from
\cite{weAiT97}, if $i\ne k$, then $G$ contains a path from $i$ to $k$, such
that the difference $(p\_{uk}-p\_{ut})$ strictly increases as $u$ progresses
from $i$ to $k$ along the path. Hence,
$p\_{kk}-p\_{kt}>p\_{ik}-p\_{it}$. Similarly,
$p\_{tt}-p\_{tk}>p\_{jt}-p\_{jk}$ whenever $j\ne t\;$. Using the above
expressions for $\D p\_{kt}$ and $\D p\_{ij}$, we get $\D p\_{kt}>\D p\_{ij}$.
\epr

\PTH{2}
Equation~(\ref{f1}) follows from the matrix-forest
theorem \cite{weAiT97} applied to the weighted multigraph $G'$ that differs
from $G$ by the weights of edges only: for all $i,j=\1n$ and $p=1\cdc
a\_{ij}$, $(\e^p_{ij})'=\tau\e^p_{ij}.$\\
\phantom{1}\hfill$\square$\par

\PPR{6}
This equality holds by virtue of the following three facts,
which are true for any $k=0\cdc n-v$ and for any $i,j,i\_1,i\_2=\1n$ such that
$i\_1\not=i\_2$:
(1) $\FF_k=\cupl^n_{i=1} \FF^{ij}_k$,
(2) $\FF^{i\_1j}_k\cap\FF^{i\_2j}_k=\varnothing$, and
(3) $\e(\FF^{ij}_k)=\e(\FF^{ji}_k) $.
\epr

\PLE{1}
Let $j\in V_i$. The desired statement follows from the
following fact: each spanning rooted forest from $\FF^{ij}_{n-v}$
can be put into correspondence with $|V_i|$ spanning rooted forests from
$\FF_{n-v}$: the latter forests have the same weight each and only differ by
the root in the component that contains $i$; each element of $\FF_{n-v}$
enters the correspondence exactly once. For $j\not\in V_i$, the statement
follows from $\FF^{ij}_{n-v}=\varnothing$.
\epr

\PPR{7}
First, we prove that $\forall\alpha\ne0,\;$ $\det(L+\alpha\oj)\ne0$.
As the matrix $L+\alpha\oj$ is reducible to a block-diagonal form, where the
blocks correspond to the connected components of $G$, it suffices to
prove its nonsingularity in the case of connected multigraphs (including the
multigraph with one vertex and without edges~--- the point graph).
Assume, on the contrary, that for some connected multigraph $G,\;$
$\det(L+\alpha\oj)=0$. Then there exists a vector ${\bf b}=(b\_1\cdc
b\_n)^{\intercal}\ne{\bf 0}$ such that $(L+\alpha\oj){\bf b}={\bf0},$ where
${\bf 0}=(0\cdc 0)^{\intercal}$. Note that the entries of $L{\bf b}$ sum to
zero, whereas the entries of $\alpha\oj{\bf b}$ are all equal. Therefore,
$L{\bf b}=\alpha\oj{\bf b}={\bf 0}$. It follows from $L{\bf b}={\bf 0}$ that
$b\_1=b\_2=\ldots=b\_n$, hence, by $\alpha\oj{\bf b}={\bf 0}$, we have
${\bf b}={\bf 0}$. This contradiction proves the invertibility of
$L+\alpha\oj$. To complete the proof, we will need a simple lemma.

\LE{2}{
For any matrices $A$ and $B$, if $A$ and $B$ are invertible and
$A\oj={\oj}B=\alpha\oj$ $(\alpha\in\R,\;\alpha\ne0),$ then
$A^{-1}\oj$ $={\oj}B^{-1}$ $={\alpha}^{-1}{\oj}.$
}

\PLE{2}
Premultiplying $A\oj=\alpha\oj$ by $A^{-1}$ yields
$\oj=\alpha A^{-1}\oj$. The statement regarding the matrix $B$ is proved
similarly.
\epr

Note that the following equalities hold true:
\begin{equation}
\label{f2}
\oj L=L\oj=0,
\end{equation}
\begin{equation}
\label{f3}
\oj^2=\oj,
\end{equation}
and, by Lemma~2 and Theorem~2, for any $\tau>0$,
\begin{equation}
\label{f4}
{(I+\tau L)}^{-1}\oj=\oj,
\end{equation}
\begin{equation}
\label{f5}
{(L+\oj)}^{-1}\oj=\oj.
\end{equation}

Using Eqs.~(\ref{f2}), (\ref{f3}), and (\ref{f5}), we obtain
\begin{equation}
\widetilde{Q}L={(L+\oj)}^{-1}L-\oj L={(L+\oj)}^{-1}(L+\oj-\oj)
=I-{(L+\oj)}^{-1}\oj=I-\oj,
\label{ImiJ}
\end{equation}
\begin{equation}
\widetilde{Q}\oj={(L+\oj)}^{-1}\oj-\oj^2=0.
\label{QtJ}
\end{equation}
Consequently, for any $\alpha\ne0,$ we have
\[
(\widetilde{Q}+\alpha^{-1}\oj)(L+\alpha\oj)=I-\oj+\oj=I,
\]
whence $\widetilde{Q}+\alpha^{-1}\oj=(L+\alpha\oj)^{-1}.$
\epr

\PPR{8}
By (\ref{ImiJ}), $\widetilde{Q}L=I-\oj$.
Similarly, $L\widetilde {Q}=I-\oj$. Thus, the first condition in the
definition of the Moore--Penrose generalized inverse is checked.
Next, using Lemma~2, (\ref{f2}), and (\ref{f3}), we have
\begin{eqnarray*}
L\widetilde{Q}L &=& L(I-\oj)=L,\\
\widetilde{Q}L\widetilde{Q} &=& (I-\oj)Q=Q-\oj Q=Q-\oj{(L+
\oj)}^{-1}+\oj^2=
Q-\oj+\oj=Q,
\end{eqnarray*}
which completes the proof.
\epr

{\bf Proof of Proposition~9} reduces to the following transformations based
on Eqs.~({\ref{f2}})--({\ref{f5}}) and Corollary~1:
\[
\quad\left(\liml_{\tau\to\infty}\tau\left({(I+\tau L)}^{-1} -
   \oj\right)+\oj\right)(L+\oj)
\]
\vspace{-1.82em}
\begin{eqnarray*}
&=&\!\!
  \liml_{\tau\to\infty}
  \tau\left({(I+\tau L)}^{-1}L+{(I+\tau L)}^{-1}\oj-\oj L -
    {\oj}^2\right)+\oj L+\oj^2\\
&=&\!\!
  \liml_{\tau\to\infty}\tau{(I+\tau L)}^{-1}L+\oj
= \liml_{\tau\to\infty}{(I+\tau L)}^{-1}(I+\tau L-I)+\oj\\
&=&\!\!
I-\liml_{\tau\to\infty}{(I+\tau L)}^{-1}+\oj=I.
\end{eqnarray*}
It now remains to apply Proposition~8.
\epr

\PTH{3}
For $j \not\in V_i$, the statement follows from Theorem~2,
Proposition~9, and the definition of $\oj$. For $j\in V_i$, using the
same and Lemma~1, we have \begin{eqnarray*} \lpij &=&
\liml_{\tau\to\infty}\tau \left({{\suml^{n-v}_{k=0}
\tau^k\e(\FF^{ij}_k)}\over
{\suml^{n-v}_{k=0}\tau^k\e(\FF_k)}}-\oj_{ij}\right)=
\liml_{\tau\to \infty} {{\suml^{n-v}_{k=0}
{\tau}^{k+1}\left(\e(\FF^{ij}_k)-{1\over |V_i|}\e(\FF_k)\right)}\over
{\suml^{n-v}_{k=0}{\tau}^k\e(\FF_k)}}\\
&=& \liml_{\tau\to \infty}
{{\suml^{n-v-1}_{k=0}{\tau}^{k+1}\left(\e(\FF^{ij}_k)-{1\over |V_i|}
\e(\FF_k)\right)}\over{\suml^{n-v}_{k=0}{\tau}^k \e(\FF_k)}}=
{{\e(\FF^{ij}_{n-v-1})-{1\over |V_i|}\e(\FF_{n-v-1})}\over
{\e(\FF_{n-v})}}.
\hbox to .208\displaywidth{\hfil\llap{$\square$}}
\end{eqnarray*}

\PPR{10}
{\em Symmetry}, {\em nonnegativity}, and
{\em disconnection condition} follow from (\ref{Qalpha}).

Let us prove {\it diagonal maximality}. The matrix $\oj$ possesses this
property in the nonstrict version $p\_{ii}\ge p\_{ij}$; therefore, by
virtue of (\ref{Qalpha}), it suffices to prove it for $Q_{n-v-1}$. By
definition, for all $i,j=\1n,\;$ $q\_{n-v-1,ij}=\e(\FF^{ij}_{n-v-1})$ holds,
where $\FF^{ij}_{n-v-1}$ is the set of all spanning rooted forests in $G$
that contain $n-v-1$ edges and have $i$ and $j$ in the same tree rooted at
$i$. Obviously, $\FF^{ij}_{n-v-1}\subseteq\FF^{ii}_{n-v-1}$.
Show that $\setminus{\FF^{ii}_{n-v-1}}{\FF^{ij}_{n-v-1}}\ne\varnothing$.
Consider an arbitrary $F\in\FF^{ij}_{n-v}$, remove from $F$ any edge that
belongs to the path from $i$ to $j$, and arbitrarily choose the root in the
newly formed component containing $j$. The resulting subgraph belongs to
$\setminus{\FF^{ii}_{n-v-1}}{\FF^{ij}_{n-v-1}}$. By the assumption of
positivity of the edge weights, we have
$\e(\FF^{ii}_{n-v-1})>\e(\FF^{ij}_{n-v-1})$, whence
$q\_{n-v-1,ii}>q\_{n-v-1,ij}$, and the property is proved. Note that
diagonal maximality can be similarly proved for $Q_1\cdc Q_{n-v-2}$;
for $Q_0$ it is obvious, whereas for $Q_{n-v}=\e(\FF_{n-v})\oj$ it is
valid in a nonstrict version.

Prove the {\it triangle inequality for proximities}. The strict
statement (for $j=k$ and $i\ne j$) follows from the diagonal maximality.
Prove that $p\_{ij}+p\_{ik}-p\_{jk}\le p\_{ii}$. For $i=j$ or $i=k$, we
have the identity. Suppose that $i\ne j$ and $i\ne k$.  Obviously,
$\FF^{ij}_{n-v-1}\cup\FF^{ik}_{n-v-1}\subseteq\FF^{ii}_{n-v-1}$.
Hence,
\begin{equation}
\e(\FF^{ij}_{n-v-1})+\e(\FF^{ik}_{n-v-1})
                    -\e(\FF^{ij}_{n-v-1}\cap\FF^{ik}_{n-v-1})
=  \e(\FF^{ij}_{n-v-1}\cup\FF^{ik}_{n-v-1})
\le\e(\FF^{ii}_{n-v-1}).
\label{simra}
\end{equation}

Define
$\FF^{ijk}_{n-v-1}$ as $\FF^{ij}_{n-v-1}\cap\FF^{ik}_{n-v-1}$ and note that
$\FF^{ijk}_{n-v-1}$ differs from
$\FF^{jik}_{n-v-1}=\FF^{ji}_{n-v-1}\cap\FF^{jk}_{n-v-1}$ only by
the roots of the trees that contain $i$, $j$, and $k$ simultaneously.
Therefore,
\begin{equation}
   \e(\FF^{ij}_{n-v-1}\cap\FF^{ik}_{n-v-1})
  =\e(\FF^{ijk}_{n-v-1})
  =\e(\FF^{jik}_{n-v-1})
\le\e(\FF^{j k}_{n-v-1}).
\label{treuh}
\end{equation}

Summing up the extreme left and extreme right parts of (\ref{simra})
and (\ref{treuh}), we obtain
\[
   \e(\FF^{ij}_{n-v-1})+\e(\FF^{ik}_{n-v-1})
\le\e(\FF^{ii}_{n-v-1})+\e(\FF^{jk}_{n-v-1}),
\]
which, by the definitions of $Q_{n-v-1}$ and $\oj$ and (\ref{Qalpha}), implies
the triangle inequality for proximities.

Prove {\it transit property}. The required inequality is valid for the
matrix $\oj$ in a nonstrict form, so by virtue of (\ref{Qalpha}), it remains
to prove it for $Q_{n-v-1}$. Obviously,
$\FF^{it}_{n-v-1}\subseteq\FF^{ik}_{n-v-1}$. To prove that
$\setminus{\FF^{ik}_{n-v-1}}{\FF^{it}_{n-v-1}}\ne\varnothing$,
consider an arbitrary $F\in\FF^{it}_{n-v}$. Remove from $F$ any edge that
belongs to the path from $k$ to $t$ and arbitrarily choose the root in the
newly formed component containing $t$. The resulting subgraph belongs to
$\setminus{\FF^{ik}_{n-v-1}}{\FF^{it}_{n-v-1}}$. By the assumption of
positivity of the edge weights, we conclude that
$\e(\FF^{ik}_{n-v-1})>\e(\FF^{it}_{n-v-1})$, and the property is proved.

To demonstrate the violation of {\it monotonicity}, it is sufficient to
consider the graph $G$ with the vertex set $V(G)=\{1,2,3\}$ and one edge
$(1,2)$ whose weight is unity. Let an edge $(1,3)$ with weight unity
be added to $G$. Here, the accessibility via dense forests provides
(for any $\alpha\ne0$)
$\D p\_{13}=-1/9<5/36=\D p\_{12}$ (which violates item~1 of
monotonicity) and $\D p\_{23}=-4/9<5/36=\D p\_{21}$ (which violates
item~2). With the same example, item~3 is also trivially violated, as
$\D p\_{22}=11/36>0$.  By adding an appropriate number of isolated
vertices, similar examples can be generated for all $n$.

\References

\end{document}